\documentclass{amsart}

\usepackage{amsmath,amssymb,amsfonts,enumerate,amsthm}

\newtheorem{thm}{Theorem}[section]
\newtheorem{cor}[thm]{Corollary}
\newtheorem{lem}[thm]{Lemma}

\newtheorem{defn}[thm]{Definition}

\newtheorem{rem}[thm]{Remark}

\begin{document}
\bibliographystyle{amsplain}
\date{}

\author{Mohammad Mahmoudi}

\address{Mohammad Mahmoudi\\Science and Research Branch, Islamic
Azad University (IAU), Tehran,
Iran.}\email{mahmoudi@damavandiau.ac.ir}

\author{Amir Mousivand}

\address{Amir Mousivand\\Science and Research Branch,Islamic Azad
University (IAU), Tehran, Iran.}\email{amirmousivand@gmail.com}

\author{Siamak Yassemi}
\address{Siamak Yassemi\\Department of Mathematics, University of
Tehran, Tehran, Iran and School of Mathematics, Institute for
Research in Fundamental Sciences (IPM), Tehran Iran.}\email{yassemi@ipm.ir}




\title{The complement of a connected bipartite graph is vertex decomposable}

\begin{abstract} Associated to a simple undirected graph $G$ is a simplicial complex $\Delta_G$ whose faces correspond
to the independent sets of $G$. A graph $G$ is called vertex
decomposable if $\Delta_G$ is a vertex decomposable simplicial
complex. We are interested in determining what families of graph
have the property that the complement of $G$, denoted
by $\overline{G}$, is vertex decomposable. We obtain the result that
the complement of a connected bipartite graph is vertex decomposable
 and so it is Cohen-Macaulay due to pureness of $\Delta_{\overline{G}}$.
\end{abstract}

\maketitle

\section{Introduction} Let $G$ be a simple graph on the vertex set $V(G)=\{ v_1 , \cdots , v_n\}$. By identifying the
vertex $v_i$ with the variable $x_i$ in the polynomial ring $k[X]=k[
x_1 , \cdots , x_n]$ over a field $k$, we can associate  to $G$ a
quadratic square-free monomial ideal $I(G)= (~ x_i x_j~|~ \{v_i ,
v_j\}\in E(G))$, where $E(G)$ is the edge set of $G$. The ideal
$I(G)$ is called the edge ideal of $G$. Using the Stanley-Reisner
correspondence, we can associate to $G$ the simplicial complex
$\Delta_G$ where $I_{\Delta_G}=I(G)$. Note that the faces of
$\Delta_G$ are the independent sets of $G$. Thus $F$ is a face of
$\Delta_G$ if and only if there is no edge of $G$ joining any two
vertices of $F$. The graph $G$ is said to be (sequentially)
Cohen-Macaulay if $k[X]/I(G)$ is a (sequentially) Cohen-Macaulay
ring.

We call a graph $G$ vertex decomposable if the simplicial complex $\Delta_G$ is vertex decomposable (see definition 2.4).
Vertex decomposability were introduced in the pure case by Provan and Billera [5] and extended to non-pure complexes by Bj\"{o}rner and Wachs [2]. We have the following implications
$$\mbox {vertex~decomposable}\Longrightarrow \mbox {shellable}\Longrightarrow \mbox {sequentially~Cohen-Macaulay}$$
and it is known that the above implications are strict.

In this article we prove that the complement (i.e. the graph whose
vertex set is $V(G)$ and edges are all the non-edges of
 $G$) of a connected bipartite graph is vertex decomposable and so shellable and sequentially Cohen-Macaulay. Since in this case $\Delta_{\overline{G}}$ is pure,
 we get the result that the complement of a connected bipartite graph is Cohen-Macaulay.

\section{Basic definitions and notations}

In this section we recall all the definitions and properties we use throughout the paper.

\begin{defn}{\bf (Complementary graph)} The Complementary graph of $G$ is the graph $\overline{G}$
with the vertex set $V(G)$ and edges all the pairs $\{ v_i , v_j\}$
such that $i\neq j$ and $\{ v_i , v_j\}\notin E$.
\end{defn}

\begin{defn}{\bf (Bipartite and complete bipartite graph)}
A bipartite graph is a graph whose vertices can be divided into two
disjoint sets $V_1$ and $V_2$ such that every edge connects a vertex
in $V_1$ to one in $V_2$. A complete bipartite graph is a bipartite
graph $G =(V_1\cup V_2,E)$ such that for any two vertices $v_1\in
V_1$ and $v_2\in V_2$, $\{v_1,v_2\}$ is an edge in G. The complete
bipartite graph with partitions of size $|V_1|=n$ and $|V_2|=m$ is
denoted by $K_{n,m}$.
\end{defn}

\begin{defn}{\bf (Cycle of graph)} A closed simple path, with no other repeated vertices than the starting and ending
vertices is called a cycle.
\end{defn}

\begin{defn} For a facet $F$ of a simplicial complex $\Delta$, the link of $F$ is the simplicial complex
$$link_{\Delta} F = \{~G~|~G\cap F =\emptyset~,~G\cup F \in \Delta~\}.$$
\end{defn}

\begin{defn}{\bf (Shedding vertex of simplicial complex and graph)} A vertex $v$ in a simplicial complex $\Delta$ is called
a shedding vertex if there is no face of $~link_{\Delta} v$ which is also a facet of $\Delta\setminus \{v\}$. A shedding vertex of a graph $G$ is the shedding vertex of the independent complex $\Delta_{G}$.
\end{defn}

\begin{defn}{\bf (Vertex decomposable simplicial complex and graph)}  A simplicial complex $\Delta$ is recursively defined to be vertex decomposable if it has only one facet or has some shedding vertex $v$ such that both $\Delta\setminus \{v\}$ and $link_{\Delta} v$ are vertex decomposable. We say that a graph $G$ is  vertex decomposable if the independent complex $\Delta_{G}$ is
vertex decomposable.
\end{defn}

\begin{rem}
Let $N_G(v)$ denotes the open neighborhood of $v$ in a graph $G$,
 i.e. all vertices adjacent to $v$, and $N_G[v]$ be the closed neighborhood of $v$ in $G$ which is
 $N_G[v]=N_G(v)\cup \{v\}$. We have the following translations of shedding vertex and vertex decomposability for the independent complex $\Delta_{G}$ (see [8, Section 2]).\\
\begin{itemize}

\item [$\bullet$] A vertex $v$ of a graph $G$ is a shedding vertex if for every independent set $S$ in $G\setminus N_G[v]$, there exists some $x\in N_G(v)$ such that $S\cup \{x\}$ is independent in $G\setminus \{v\}$.\\

\item [$\bullet$] A graph $G$ is vertex decomposable if it is a discrete graph or has some shedding vertex $v$ such that both $G\setminus \{v\}$ and $G\setminus N_G[v]$ are vertex decomposable.

\end{itemize}
\end{rem}

\begin{rem} Recall that a vertex $v$ in a graph $G$ is called simplicial vertex if $N_G[v]$ is a clique of $G$. In [8], Woodroofe  showed that any neighbor of a simplicial vertex is a shedding vertex for $G$ and that any chordal graph is vertex decomposable. Therefore any complete graph is  vertex decomposable.
\end{rem}

\section{Main result}
In this section we state and prove the main theorem of this paper that says the complement of
a connected bipartite graph is vertex decomposable. We split the proof into some special cases. First we prove the result in the case where $G$ has a free vertex (vertex of degree 1). Then we focus on the problem with the assumption that $G$ has no free vertex and conclude that it would contain a shedding vertex. Finally we use the fact that $G$ contains at least a shedding vertex, say $x$, with the property that $\overline{G}\setminus \{x\}$ is the complement of a connected bipartite
graph and apply the induction.

\begin{lem} Let $G$ be a connected bipartite graph. Suppose $V(G)=V_1\cup V_2$ and there exists $v\in V(G)$ such
that $N_G(v)=V_1$ or $N_G(v)=V_2$. Then $\overline{G}$ is vertex decomposable.
\end{lem}

\begin{proof} let $V_1=\{ x_1 , \cdots , x_n\}$ and $V_2=\{ y_1 , \cdots , y_m\}$. If $n=1$ or $m=1$,
then the connected components of $\overline{G}$ are  vertex
decomposable. Therefore $\overline{G}$ is vertex decomposable, cf.
[8, Lemma 6.1]. Assume  that $n\geq 2$ and $m\geq 2$. We may assume
$N_G(x_1)=V_2$. Therefore $N_{\overline{G}}[x_1]=V_1$ and hence it
is a clique of $\overline{G}$. This means that $x_1$ is a simplicial
vertex of $\overline{G}$. Therefore $\overline{G}$ has a shedding
vertex, say $x_j$, with $j>1$. We claim $\overline{G}\setminus
\{x_j\}$ and $\overline{G}\setminus N_{\overline{G}}[x_j]$ are
vertex decomposable. We have that $\overline{G}\setminus
\{x_j\}=\overline{G\setminus \{x_j\}}$ and $G\setminus \{x_j\}$ is a
bipartite graph which is also connected, since $N_G(x_1)=V_2$. Thus
by induction hypothesis $\overline{G\setminus \{x_j\}}$ is vertex
decomposable. On the other hand, $\overline{G}\setminus
N_{\overline{G}}[x_j]$ is a complete graph over a subset of $V_2$,
and so is vertex decomposable.
\end{proof}

\begin{cor} The complement of the graph $K_{n,m}$ is vertex decomposable.
\end{cor}

\begin{lem} Let $G$ be a connected bipartite graph without free vertex (so it would contains a cycle). Then there exist some cycle $C$ of $G$ and some $v\in V(C)$ such that $G\setminus\{v\}$ is connected.
\end{lem}

\begin{proof} Suppose the contrary that for each vertex $v$ of any cycle of $G$, the bipartite graph
$G\setminus\{v\}$ is disconnected. Let the number of cycles in $G$
is $t$ and let $\{x_{11} , x_{12} , \cdots , x_{1n_1}\}$, $\{x_{21}
, x_{22} , \cdots , x_{2n_2}\}$, $\cdots$ ,
$\{x_{t1} , x_{t2} , \cdots , x_{tn_t}\}$  be all the cycles of $G$.\\
$G\setminus\{x_{1n_1}\}$
is disconnected, so suppose
$$G\setminus\{x_{1n_1}\}=G_{11}\cup G_{12}\cup\cdots\cup G_{1l_1}$$
be the decomposition of $G\setminus\{x_{1n_1}\}$ as the union of its connected components.
We may assume $\{x_{11} , x_{12} , \cdots , x_{1(n_1-1)}\}\subseteq V(G_{11})$.
Since $G$ is connected, there exists $\alpha\in V(G_{12})\cup\cdots\cup V(G_{1l_1})$ such that $\alpha$ is adjacent to
$x_{1n_1}$ and we may assume $\alpha\in V(G_{12})$.
We have that $\mbox{deg}_G(\alpha)>1$, and so there exists $\beta\in V(G)\setminus \{x_{1n_1}\}$ such that $\beta$ is
adjacent to $\alpha$. It is easy to see that $\beta\notin V(G_{11})$. Again we have that $\mbox{deg}_G(\beta)>1$.
Proceeding in this way (bringing in the mind that $G$ has no free vertex), we would obtain a cycle which has no
intersection with $V(G_{11})$. Let $\{x_{21} , x_{22} , \cdots , x_{2n_2}\}$ be the described cycle.\\
Similarly, suppose $$G\setminus\{x_{2n_2}\}=G_{21}\cup G_{22}\cup\cdots\cup G_{2l_2}$$ be the decomposition of
$G\setminus\{x_{2n_2}\}$ as the union of its connected components such that
$\{x_{21} , x_{22} , \cdots , x_{2(n_2-1)}\}\subseteq V(G_{21})$. Since $G$ is connected, there exists
$\alpha'\notin V(G_{21})$ such that $\alpha'$ is adjacent to $x_{2n_2}$, but $\mbox{deg}_G(\alpha')>1$ and so
there exists $\beta'\notin V(G_{21})$ which is adjacent to $\alpha'$. Proceeding in this way provides a cycle
which has no intersection with $V(G_{11})\cup V(G_{21})$.\\
If we continue the above described procedure, after $t-1$ stage, we get that the cycle
$\{x_{t1} , x_{t2} , \cdots , x_{tn_t}\}$ has no intersection with $V(G_{11})\cup V(G_{21})\cup\cdots\cup V(G_{(t-1)1})$.\\
Let $$G\setminus\{x_{tn_t}\}=G_{t1}\cup G_{t2}\cup\cdots\cup G_{tl_t}$$ be the decomposition of $G\setminus\{x_{tn_t}\}$
as the union of its connected components such that $\{x_{t1} , x_{t2} , \cdots , x_{t(n_t-1)}\}\subseteq V(G_{t1})$.\\
A similar argument as above shows that there exists $\alpha''\notin V(G_{t1})$ which is adjacent to $x_{tn_t}$.
Now $\mbox{deg}_G(\alpha'')>1$ implies that there exists a cycle in $G$ that has no intersection with
$V(G_{11})\cup V(G_{21})\cup\cdots\cup V(G_{t1})$ which is impossible.
\end{proof}

The strategy of the proof of our main theorem is based on the existence of free vertex. In the case where $G$ does not
contain free vertex, we will obtain the result that any vertex of each cycle of $G$ is a shedding vertex, and moreover,
it follows from previous lemma that there exists at least one vertex in a cycle of $G$, say $x$, such that
$G\setminus \{x\}$ is connected. Finally we use induction to conclude the result.

\begin{thm} 
The complement of a connected bipartite graph is vertex decomposable and so it is Cohen-Macaulay.
\end{thm}

\begin{proof} Let $G$ be a connected bipartite graph and suppose that $V(G)=V_1\cup V_2$ where
$V_1=\{ x_1 , \cdots , x_n\}$ and $V_2=\{ y_1 , \cdots , y_m\}$.
 In view of Lemma 3.1, we may assume  $n\geq 2$, $m\geq 2$ and that $N_G(x_i)\neq V_2$ and $N_G(y_j)\neq V_1$ for all
$i=1 , \cdots , n$ and $j= 1 , \cdots , m$. We split the argument into two cases.\\
{\bf Case (1)} Assume that $G$ has a free vertex, say $x_1$.\\
First we show that $x_1$ is a shedding vertex of $\overline{G}$. Let
$S$ be an independent subset of $\overline{G}\setminus
N_{\overline{G}}[x_1]$. We have to find $v\in N_{\overline{G}}(x_1)$
such that $S\cup \{v\}$ is an independent subset of
$\overline{G}\setminus \{x_1\}$. If $S=\emptyset$, then it follows
from $N_G(x_1)\neq V_2$ that $V_2\cap N_{\overline{G}} (x_1)\neq
\emptyset$, and there is nothing to prove.\\
 Assume $S\neq
\emptyset$. We show that $|S|=1$. Suppose the contrary that $|S|>1$
and let $u,w\in S$. Therefore $\{u,w\}$ is an edge of $G$ and hence
we may assume $u\in V_1$ and $w\in V_2$. This implies that
$u\in N_{\overline{G}}[x_1]$ which is contradiction. Therefore $|S|=1$. We know that
$V_1\subseteq N_{\overline{G}}[x_1]$, hence $S\subseteq V_2$ and thus we may assume $S=\{y_1\}$.\\
We claim that $N_G(y_1)\setminus\{x_1\}\neq\emptyset$. Suppose in contrary that $N_G(y_1)=\{x_1\}$.
It follows that $\mbox{deg}(y_1)=1$ which together with $\mbox{deg}(x_1)=1$, $n\geq 2$, and $m\geq 2$
implies that the edge $\{x_1,y_1\}$ is a connected component of $G$ which is impossible.\\
Let $x_2\in N_G(y_1)\setminus\{x_1\}$. Then $\{x_2,y_1\}$ is not an
edge of $\overline{G}$ and $x_2\in N_{\overline{G}} (x_1)$. So
$S\cup\{x_2\}=\{x_2,y_1\}$ is an independent subset of
$\overline{G}$ and hence $x_1$ is a shedding vertex of
$\overline{G}$. The next step is to show that
$\overline{G}\setminus\{x_1\}$ and $\overline{G}\setminus
N_{\overline{G}}[x_1]$ are vertex decomposable. $G\setminus\{x_1\}$
is a bipartite graph which is also connected because
$\mbox{deg}(x_1)=1$. Thus $\overline{G\setminus \{x_1\}}$ is vertex
decomposeble by induction hypothesis. Vertex decomposability of
$\overline{G}\setminus\{x_1\}$ follows from
$\overline{G}\setminus\{x_1\}=\overline{G\setminus \{x_1\}}$. Note
that $V_1\subseteq N_{\overline{G}}[x_1]$ and so
$\overline{G}\setminus N_{\overline{G}}[x_1]$ is a complete graph
over
a subset of $V_2$ which is vertex decomposable.\\
{\bf Case (2)} Assume that $G$ has no free vertex.\\
In this case $G$ contains at least a cycle. Suppose $x_1$ belongs to
a cycle of $G$. We claim that $x_1$ is a shedding vertex of
$\overline{G}$. Let $S$ be an independent subset of
$\overline{G}\setminus N_{\overline{G}}[x_1]$. A similar argument as
in the Case (1) implies that $|S|=1$ and hence we may assume that
$S=\{y_1\}$. Since $G$ has no free vertex, we get that
$N_G(y_1)\setminus\{x_1\}\neq\emptyset$. Let $x_2\in
N_G(y_1)\setminus\{x_1\}$. Then $x_2\in N_{\overline{G}} (x_1)$ and
$S\cup\{x_2\}=\{x_2,y_1\}$ is an independent subset of
$\overline{G}$. Therefore $x_1$ is a shedding vertex of
$\overline{G}$. Note that this argument shows that any vertex of $G$
appeared in a cycle is a shedding vertex of $\overline{G}$. Now
suppose that $x_1$ be as in Lemma 3.3. To complete the proof, it is
enough to show that $\overline{G}\setminus\{x_1\}$ and
$\overline{G}\setminus N_{\overline{G}}[x_1]$ are vertex
decomposable. The result for $\overline{G}\setminus
N_{\overline{G}}[x_1]$ is similar to the Case (1). Finally it
follows from Lemma 3.3 that $G\setminus\{x_1\}$ is a connected
bipartite graph and hence
$\overline{G}\setminus\{x_1\}=\overline{G\setminus \{x_1\}}$ is
vertex decomposable by induction hypothesis.
\end{proof}

\begin{cor} The complement of any cycle of even length is vertex
decomposable.
\end{cor}

\begin{center}
\end{center}

\end{document}